\documentclass[12pt]{article}
\usepackage{amssymb}
\usepackage{amsmath,amsthm}
\usepackage[latin1]{inputenc}
\usepackage{citesort}
\usepackage{graphicx}
 \newtheorem{remark}{Remark}

 \newtheorem{lemma}[remark]{Lemma}
 \newtheorem{theorem}[remark]{Theorem}
 
 \newtheorem{corollary}[remark]{Corollary}

\title{On defensive alliances and line graphs}

\author{J. M. Sigarreta\footnote{e-mail:\mbox{\tt
    josemaria.sigarreta\@@uc3m.es}}\\
{\em Department of Mathematics }\\ Carlos III  University of Madrid\\
Avda. de la Universidad 30, 28911
Leganés (Madrid),  Spain \\
J. A. Rodr\'{\i}guez\footnote{e-mail:\mbox{\tt
    juanalberto.rodriguez\@@urv.net}} \\
{\em Department of Computer Engineering and Mathematics}\\
Rovira i Virgili University of Tarragona\\ Av. Pa\"{\i}sos Catalans
26, 43007 Tarragona, Spain}

\date{}

\begin{document}

\maketitle

\begin{abstract}
Let  $\Gamma$ be a simple graph of size $m$ and degree sequence
$\delta_1\ge \delta_2\ge \cdots \ge \delta_n$. Let ${\cal
L}(\Gamma)$ denotes the line graph of $\Gamma$. The aim of this
paper is to study mathematical properties of the alliance number,
${a}({\cal L}(\Gamma)$, and the global alliance number,
$\gamma_{a}({\cal L}(\Gamma))$, of the line graph of a simple graph.
We show that
$\left\lceil\frac{\delta_{n}+\delta_{n-1}-1}{2}\right\rceil \le
{a}({\cal L}(\Gamma))\le \delta_1.$ In particular,
 if $\Gamma$ is a $\delta$-regular graph ($\delta>0$),
then $a({\cal L}(\Gamma))=\delta$, and if $\Gamma$ is a
$(\delta_1,\delta_2)$-semiregular bipartite graph, then $a({\cal
L}(\Gamma))=\left\lceil \frac{\delta_1+\delta_2-1}{2} \right\rceil$.
As a consequence of the study we compare $a({\cal L}(\Gamma))$ and
${a}(\Gamma)$, and we characterize the graphs having $a({\cal
L}(\Gamma))<4$. Moreover, we show that the global-connected
allian\-ce number of ${\cal L}(\Gamma)$ is bounded by $\gamma_{ca}
\left( {\cal L}(\Gamma)\right) \ge
\left\lceil\sqrt{D(\Gamma)+m-1}-1\right\rceil,$ where $D(\Gamma)$
denotes the diameter of $\Gamma$, and we show that the global
alliance number of ${\cal L}(\Gamma)$ is bounded by
$\gamma_{a}({\cal L}(\Gamma))\geq
\left\lceil\frac{2m}{\delta_{1}+\delta_{2}+1}\right\rceil$. The case
of strong alliances is studied by analogy.

\end{abstract}

{\it Keywords:}  Defensive alliance,  alliances in graphs, line
graph.

{\it AMS Subject Classification numbers:}   05C69;  15C05

\section{Introduction}

The study of defensive alliances in graphs, together with a variety
of other kinds of allian\-ces,  was introduced in
\cite{alliancesOne}. In the referred paper was initiated the study
of the mathe\-matical properties of alliances. In particular,
several bounds on the defensive alliance number were given. The
particular case of global (strong) defensive alliance was
investigated in\cite{GlobalalliancesOne} where several bounds on the
global (strong) defensive alliance number were obtained.

In  \cite{spectral} were obtained several tight bounds on different
types of allian\-ce numbers of a graph, namely (global) defensive
alliance number, (global) offensive alliance number and (global)
dual alliance number. In particular, was investigated the
relationship between the alliance numbers of a graph and its
algebraic connectivity, its spectral radius, and its Laplacian
spectral radius. A particular study of the alliance numbers, for the
case of planar graphs, can be found in \cite{planar}.  Moreover, for
the study of offensive alliances we cite \cite{favaron,offensive}.

The aim of this paper is to study mathematical properties of the
alliance number and the global alliance number of the line graph of
a simple graph. We begin by stating some notation and terminology.
In this paper $\Gamma=(V,E)$ denotes a simple graph of order $n$ and
size $m$. The degree sequence of $\Gamma$ will be denoted by
$\delta_1\ge \delta_2\ge \cdots \ge \delta_n$. Moreover, the degree
of a vertex $v\in V$ will be denoted by $\delta(v)$. The line graph
of $\Gamma$ will be denoted by ${\cal L}(\Gamma)=(V_l,E_l)$. The
degree of the vertex $e=\{u,v\}\in V_l$ is
$\delta(e)=\delta(u)+\delta(v)-2$. The subgraph induced by a set
$S\subset V$ will be denoted by $\langle S\rangle$.

For a non-empty subset $S\subseteq V$, and any vertex $v\in V$, we
denote by $N_S(v)$ the set of neighbors $v$ has in $S$:
$N_S(v):=\{u\in S: u\sim v\},$
 Similarly, we denote by
$N_{V\setminus S}(v)$ the set of neighbors $v$ has in $V\setminus
S$: $N_{V\setminus S}(v):=\{u\in V\setminus S: u\sim v\}$.

A nonempty set of vertices $S\subseteq V$ is called a {\em defensive
allian\-ce} if for every $v\in S$, $| N_S(v) | +1\ge  |
N_{V\setminus S}(v)|.$ Equivalently, $S$ is a defensive alliance if
for every $v\in S$, $2| N_S(v) | +1\ge \delta(v)$. In this case, by
strength of numbers, every vertex in $S$ is {\em defended} from
possible attack by vertices in $V\setminus S$. A defensive alliance
$S$ is called {\em strong} if for every $v\in S$, $| N_S(v) | \ge |
N_{V\setminus S}(v)|$. Equivalently, $S$ is a defensive alliance if
for every $v\in S$, $2| N_S(v) | \ge \delta(v)$. In this case every
vertex in $S$ is {\em strongly defended}.

The {\em defensive alliance number} $a(\Gamma)$ (respectively, {\em
strong defensive allian\-ce number} $\hat{a}(\Gamma)$) is the
minimum cardinality of any defensive alliance (respectively, strong
defensive alliance) in $\Gamma$. A defensive alliance, $S$, in
$\Gamma$ is \emph{minimal} if no proper subset of $S$ is a defensive
alliance. A \emph{minimum} defensive alliance is a minimal defensive
alliance  of smallest cardinality, i.e., $|S|=a(\Gamma)$.

A particular case of alliance, called global defensive allian\-ce,
was studied in \cite{GlobalalliancesOne}. A defensive alliance $S$
is called {\em global} if it affects every vertex in $V\setminus S$,
that is, every vertex in $V\setminus S$ is adjacent to at least one
member of the alliance $S$. Note that, in this case, $S$ is a
dominating set.  The {\em global defensive allian\-ce number}
$\gamma_a(\Gamma)$ (respectively, {\em global strong defensive
alliance number} $\gamma_{\hat{a}}(\Gamma)$) is the minimum
cardinality of any global defensive alliance (respectively, global
strong defensive alliance) in $\Gamma$. Singular interest displays
the global defensive alliances whose induced subgraph is connected.
We define the {\em global-connected defensive allian\-ce number},
 $\gamma_{ca}(\Gamma)$,
 (respectively, {\em global-connected strong defensive
alliance number} $\gamma_{c\hat{a}}(\Gamma)$) as the minimum
cardinality of any global defensive alliance (respectively, global
strong defensive alliance) in $\Gamma$ whose induced subgraph is
connected.

In this paper we show that the alliance number of ${\cal L}(\Gamma)$
is bounded by
$\left\lceil\frac{\delta_{n}+\delta_{n-1}-1}{2}\right\rceil \le
{a}({\cal L}(\Gamma))\le \delta_1$. In particular,
 if $\Gamma$ is a $\delta$-regular graph ($\delta>0$),
then $a({\cal L}(\Gamma))=\delta$, and if $\Gamma$ is a
$(\delta_1,\delta_2)$-semiregular bipartite graph, then $a({\cal
L}(\Gamma))=\left\lceil \frac{\delta_1+\delta_2-1}{2} \right\rceil$.
As a consequence of the study we compare $a({\cal L}(\Gamma))$ and
${a}(\Gamma)$, and we characterize the graphs having $a({\cal
L}(\Gamma))<4$. In the case of global alliances, we show that the
global alliance number of ${\cal L}(\Gamma)$  is bounded by
$\gamma_{a}({\cal L}(\Gamma))\geq
\left\lceil\frac{2m}{\delta_{1}+\delta_{2}+1}\right\rceil$ and the
global-connected alliance number of ${\cal L}(\Gamma)$ is bounded by
$\gamma_{ca} \left( {\cal L}(\Gamma)\right) \ge
\left\lceil\sqrt{D(\Gamma)+m-1}-1\right\rceil,$ where $D(\Gamma)$
denotes the diameter of $\Gamma$. In addition, the case of strong
alliances is studied by analogy.

\section{Defensive alliances  and line graphs}

\begin{theorem}\label{thcota}
Let $\Gamma$ be a graph whose degree sequence is $\delta_{1}\geq
\delta_{2}\geq \cdots \geq\delta_{n}$. Then
%\begin{equation}
$\left\lceil\frac{\delta_{n}+\delta_{n-1}}{2}\right\rceil \le
\hat{a}({\cal L}(\Gamma))\le \delta_1$ and
%\end{equation}
%\begin{equation}\label{minimal}
$\left\lceil\frac{\delta_{n}+\delta_{n-1}-1}{2}\right\rceil \le
a({\cal L}(\Gamma))\le \delta_1.$
%\end{equation}
 Moreover, if  $\Gamma$ has a unique vertex of maximum degree, then
%\begin{equation}
$a({\cal L}(\Gamma))\le  \delta_1-1.$
%\end{equation}
\end{theorem}

\begin{proof}

If  $S_l$ denotes a  strong defensive alliance in ${\cal
L}(\Gamma)$, then   $\forall e\in S$, $$2(|S_l|-1)\ge
2|N_{S_l}(e)|\ge \delta(e)\ge \delta_{n}+\delta_{n-1}-2.$$
Therefore, the lower bound of $\hat{a}({\cal L}(\Gamma))$ follows.

Let $v$ be a vertex of maximum degree in $\Gamma$ and let
$S_v=\{e\in E: v\in e\}$. Thus, $\langle S_v\rangle \cong
K_{\delta_1}$ and, as a consequence,  $\forall e\in S_v$,
$|N_{S_v}(e)|=\delta_1-1\ge \delta_2-1\ge |N_{V_l\backslash
S_v}(e)|$. Hence, $S_v\subset V_l$ is a strong defensive alliance in
${\cal L}(\Gamma)$.  So, $\hat{a}({\cal L}(\Gamma))\le \delta_1.$

 The lower bound of ${a}({\cal L}(\Gamma))$ is obtained by analogy to the previous case.
 Moreover,  ${a}({\cal L}(\Gamma))\le
\hat{a}({\cal L}(\Gamma))\le \delta_1$.

Suppose that $v\in V$ is the  unique vertex of maximum degree in
$\Gamma$. As above, let $S_v=\{e\in E: v\in e\}$. Let $e'\in S_v$
and let  $S_v'=S_v\backslash \{e'\}$. Thus, $\langle S_v'\rangle
\cong K_{\delta_1-1}$. Hence,  $\forall e\in S_v'$,
$|N_{S_v'}(e)|+1=\delta_1-1\ge \delta_2\ge |N_{V_l\backslash
S_v'}(e)|$. Therefore, $S_v'\subset V_l$ is a defensive alliance in
${\cal L}(\Gamma)$ and its cardinality is $\delta_1-1$. So, $a({\cal
L}(\Gamma))\le \delta_1-1.$
\end{proof}

\begin{corollary}
If  $\Gamma$ is a $\delta$-regular graph ($\delta>0$), then $a({\cal
L}(\Gamma))=\hat{a}({\cal L}(\Gamma))=\delta.$
\end{corollary}

%\begin{proof}
%Sea $\Gamma=(V,E)$ y sea $S_v=\{a\in E: v\in a\}$. Entonces
%$|N_{S_v}(a)|=\delta-1$ y $|N_{V_l\backslash S_v}(a)|=\delta-1$.
%Así, $S_v\subset V_l$ es una alianza defensiva fuerte en ${\cal
%L}(\Gamma)$. Supongamos que existe una alianza defensiva $S\subset
%E$ tal que $|S|< \delta$. Si $a\in S$, entonces $|N_{S}(a)|\le
%\delta-2$. Por tanto, $|N_{V_l\backslash S}(a)|\ge \delta$. Lo que
%contradice lo supuesto.
%\end{proof}

\begin{theorem}
If $\Gamma$ is a  $(\delta_1,\delta_2)$-semiregular bipartite graph,
then
$$a({\cal L}(\Gamma))=\left\lceil \frac{\delta_1+\delta_2-1}{2}
\right\rceil \quad {\rm and} \quad \hat{a}({\cal
L}(\Gamma))=\left\lceil\frac{\delta_1+\delta_2}{2}\right\rceil .$$
\end{theorem}

\begin{proof}
Suppose $\delta_1>\delta_2$. By Theorem \ref{thcota}, we only need
to show that there exists a defensive alliance whose cardinality is
$\left\lceil \frac{\delta_1+\delta_2-1}{2} \right\rceil$. Let $v\in
V$ be a vertex of maximum degree in $\Gamma$ and let $S_v=\{e\in E:
v\in e\}$. Hence, $\langle S_v \rangle \cong K_{\delta_1}$.
Therefore, taking $S\subset S_v$ such that $|S|=\left\lceil
\frac{\delta_1+\delta_2-1}{2} \right\rceil$, we obtain $\langle
S\rangle \cong K_{\left\lceil \frac{\delta_1+\delta_2-1}{2}
\right\rceil}$. Thus,  $\forall e\in S$, $$|N_{S}(e)|+1=\left\lceil
\frac{\delta_1+\delta_2-1}{2} \right\rceil\ge
\delta_1+\delta_2-1-\left\lceil \frac{\delta_1+\delta_2-1}{2}
\right\rceil = |N_{V_l\backslash S}(e)|.$$ So, $S$ is a defensive
alliance in  ${\cal L}(\Gamma)$. The proof of $\hat{a}({\cal
L}(\Gamma))=\left\lceil\frac{\delta_1+\delta_2}{2}\right\rceil $ is
analogous to the previous one.
\end{proof}

%\begin{proof}
%Sea $S_l\subset V_l$ una alianza defensiva de cardinal mínimo en
%${\cal L}(\Gamma)$. Probaremos que el conjunto  $S\subset V$
%definido como
%$$S=\{v\in V: v\in a, \mbox{ } \mbox{\rm para algún }\mbox{ } a\in
%S_l\}$$ contiene una alianza defensiva $S_1$ tal que
%$|S_1|\le|S_l|.$

%Nótese que como la alianza $S_l$ es minimal, entonces $\langle
%S_l\rangle$ es conexo y, por consiguiente, $\langle S\rangle$
%también es conexo. Por tanto,  $|S|=|S_l|+1$.

%La primera parte, ya está.

%Sea  $a=\{u,v\}\in S_l$. Como $S_l$ es una alianza defensiva en
%${\cal L}(\Gamma)$, se cumple
%\begin{equation}\label{Si}
%2|N_{S_l}(a)|+1\geq \delta(a).
%\end{equation}
%Por otro lado, $|N_{S_l}(a)|=|N_{S}(u)|+|N_{S}(v)|-2$ y $
% \delta (a)=\delta (u)+\delta(v)-2$. Por lo tanto,
%\begin{equation}\label{Su}
%2|N_{S}(u)|+2|N_{S}(v)|-1 \geq \delta (u)+\delta (v).
%\end{equation}

%Nótese que:
%\begin{equation}\label{Igual}
%|N_{S_{1}}(u)|= |N_{S_{1}\backslash {v}}(u)|+1,
%\end{equation}
%y que:
%\begin{equation}\label{Igual1}
%|N_{S_{1}}(v)|= |N_{S_{1}\backslash {u}}(v)|+1.
%\end{equation}
%Luego de (\ref{Ssss}) y (\ref{Igual}) se tiene que:

%\begin{equation}\label{bien}
%2|N_{S_{1}}(u)|\leq \delta(u).
%\end{equation}

Now we are going to characterize the graphs having $a({\cal
L}(\Gamma))< 4$.

\begin{lemma}{\rm \cite{alliancesOne}}  \label{lemmaCaract}
For any graph $\Gamma$,
\begin{enumerate}
\item $a(\Gamma)=1$ if and only if there exists a vertex $v\in V$
such that $\delta(v)\le 1$.

\item $a(\Gamma)=2$ if and only if  $2\le \displaystyle\min_{v\in V}\{\delta(v)\}$ and  $\Gamma$
has two adjacent vertices of degree at most three.

\item $a(\Gamma)=3$ if and only if  $a(\Gamma)\neq 1$, $a(\Gamma)\neq
2$, and  $\Gamma$ has an induced subgraph isomorphic to either (a)
  $P_{3}$, with vertices, in order, $u$, $v$ and $w$, where
$\delta(u)$ and  $\delta(w)$ are at most three, and $\delta(v)$ is
at most five, or (b) isomorphic to $K_{3}$, each vertex of which has
degree at most five.
\end{enumerate}

\end{lemma}

\begin{theorem}
For any graph $\Gamma$,

\begin{enumerate}
\item  $a({\cal L}(\Gamma))=1$ if and only if either $\Gamma$ has a
connected component isomorphic to  $K_2$, or $\Gamma$ has a vertex
of degree one which is adjacent to a vertex of degree two.

\item   $a({\cal L}(\Gamma))=2$ if and only if   $a({\cal L}(\Gamma))\ne
1$ and  $\Gamma$ has a subgraph isomorphic to $P_3$, with vertices,
in order, $u,v$ and $w$, such that $\delta(u)+\delta(v)\le 5$ and
$\delta(v)+\delta(w)\le 5$.

\item $a({\cal L}(\Gamma))=3$ if and only if $a({\cal L}(\Gamma))\ne 1$,  $a({\cal L}(\Gamma))\ne 2$,
and $\Gamma$ has a subgraph isomorphic to either  (a) $P_4$, with
vertices, in order, $u,v,w$ and $x$, such that
$\delta(u)+\delta(v)\le 5$, $\delta(x)+\delta(w)\le 5$ and
$\delta(v)+\delta(w)\le 7$, or  (b) $K_3$, with vertices
$\{u,v,w\}$, such that
 $\delta(u)+\delta(v)\le 7$,  $\delta(u)+\delta(w)\le 7$ and
$\delta(v)+\delta(w)\le 7$, or (c) $K_{1,3}$, with vertices
$\{u,v,w,x\}$, and hub $v$, such that
 $\delta(v)+\delta(u)\le 7$,  $\delta(v)+\delta(w)\le 7$  and
$\delta(v)+\delta(x)\le 7$.
\end{enumerate}
\end{theorem}

\begin{proof}

 The result follows from Lemma \ref{lemmaCaract}:

\begin{enumerate}
\item   ${\cal L}(\Gamma)$ has an isolated vertex if and only if $\Gamma$ has a
connected component isomorphic to  $K_2$.  Moreover, ${\cal
L}(\Gamma)$ has a vertex of degree one if and only if $\Gamma$ has a
vertex of degree one adjacent to a vertex of degree two.

\item  ${\cal L}(\Gamma)$ has two adjacent vertices,  $e_1,e_2 \in
V_l$, such that  $\delta(e_1)\le 3$ and  $\delta(e_2)\le 3$, if and
only if $\Gamma$ has three vertices $u,v,w\in V$ such that
$e_1=\{u,v\}$ and $e_2=\{v,w\}$, with
$\delta(u)+\delta(v)-2=\delta(e_1)\le 3$ and
$\delta(v)+\delta(w)-2=\delta(e_2)\le 3$.

\item ${\cal L}(\Gamma)$ has an induced subgraph isomorphic to $P_{3}$, with vertices, in order, $e_1$, $e_2$ and
 $e_3$, where $\delta(e_1)\le 3$,
$\delta(e_2)\le 3$ and $\delta(e_3)\le 5$  if and only if $\Gamma$
has a subgraph  isomorphic to $P_4$, with vertices, in order,
$u,v,w$ and $x$, where  $e_1=\{u,v\}$, $e_2=\{v,w\}$, $e_3=\{w,x\}$,
$\delta(u)+\delta(v)\le 5$, $\delta(x)+\delta(w)\le 5$ and
$\delta(v)+\delta(w)\le 7$.

On the other hand, ${\cal L}(\Gamma)$ has an induced subgraph
isomorphic to  $K_{3}$ if and only if either $\Gamma$ has a subgraph
isomorphic to $K_{3}$, or $\Gamma$ has a subgraph isomorphic to
$K_{1,3}$. Moreover, for $e=\{u,v\}\in V_l$, $\delta(e)\le 5$ if and
only if $\delta(u)+\delta(v)\le 7$.
\end{enumerate}

\end{proof}

We remark that a similar characterization can be done in the case of
strong alliances.

Now we are going to compare  $a(\Gamma)$ and $ a({\cal L}(\Gamma))$.
There are cases in which  $a(\Gamma)= a({\cal L}(\Gamma)).$ A
trivial instance is the case $\Gamma\cong C_k$ ($\Gamma$ isomorphic
to the cycle of length $k$). In order to show the case $a({\cal
L}(\Gamma))< a(\Gamma)$ we take $\Gamma\cong O_5$ (the odd graph
$O_5$). That is, $a({\cal L}(O_5)=5<
6=girth(O_5)=a(O_5)$\footnote{It was shown in \cite{alliancesOne}
that if $\Gamma$ is 5-regular, then ${a}(\Gamma)=girth(\Gamma)$.}.
Moreover, there are cases in which $a(\Gamma)< a({\cal L}(\Gamma))$.
For instance, if either $\Gamma$ is isomorphic to a tree, or
$\Gamma$ is isomorphic to an unicyclic\footnote{A connected graph
containing exactly one cycle.} graph, but $\Gamma\ncong C_k$, then
$1=a(\Gamma)\le a({\cal L}(\Gamma))$. In particular, if $\Gamma\cong
K_{1,n}$, $n>2$, then $a(\Gamma)=1< n-1= a({\cal L}(\Gamma))$.

\begin{figure}[h]
\begin{center}
\caption{$\Gamma$ and its line graph ${\cal L}(\Gamma)$ }
\label{figure}
\vspace{-3cm}
\includegraphics[angle=0, width=5.5cm]{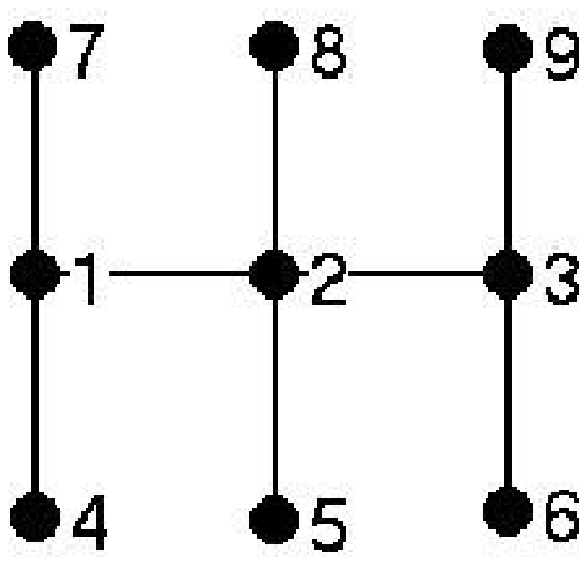}
\includegraphics[angle=0, width=5.5cm]{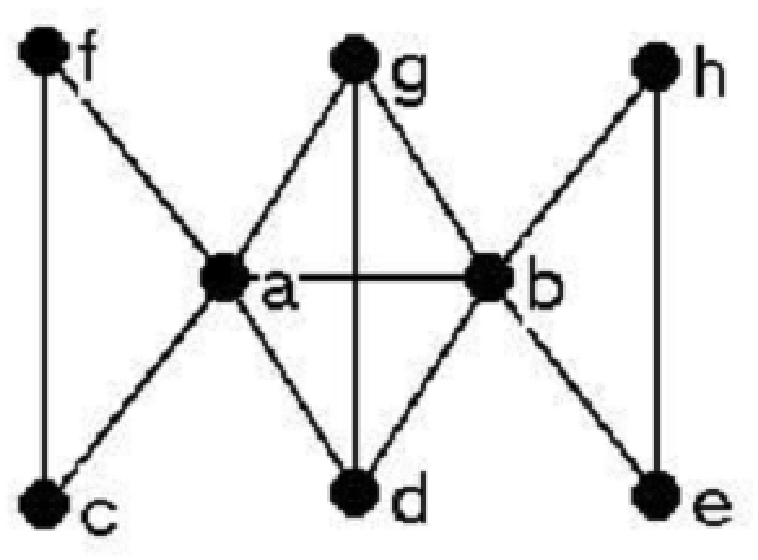}
\end{center}
\end{figure}
\vspace*{-3cm}

%\begin{figure}[h]
%\begin{center}
%\includegraphics[width=0.2\textwidth]{G}
%\hspace{0,5cm}
%\includegraphics[width=0.25\textwidth]{LG}
%\end{center}
%\end{figure}

We define the \emph{characteristic set\/} of    $S_l \subset V_l$ as
$C_{S_l}:=\{v\in V: v\in e, \mbox{ \rm for some } e\in S_l\}.$ For
instance, in the graph of Figure \ref{figure}, $S_l=\{f,c\}$ is a
minimum  defensive alliance in ${\cal L}(\Gamma)$ and its
characteristic set, $C_{S_l}=\{1,4,7\}$, is a  defensive alliance in
$\Gamma$. Notice that $C_{S_l}$ contains the  defensive alliances
$S_1=\{1,4\}$, $S_2=\{1,7\}$, $S_3=\{4\}$ and $S_4=\{7\}$.  We
emphasize that in some cases ${\cal L}(\Gamma)$ has not minimum
defensive alliances such that its characteristic set is a defensive
alliance in $\Gamma$.

\begin{theorem}
If there exists a minimum defensive alliance in ${\cal L}(\Gamma)$
such that its characteristic set is a defensive alliance in
$\Gamma$, then $a(\Gamma)\leq a({\cal L}(\Gamma)).$
\end{theorem}

\begin{proof}

Let $S_l\subset V_l$ be a minimum defensive alliance in ${\cal
L}(\Gamma)$.  We shall show that the characteristic set of $S_l$,
$C_{S_l}$, contains a defensive alliance whose cardinality is $\le
|S_l|.$

As $S_{l}$ is a minimum defensive alliance in ${\cal L}(\Gamma)$,
then the subgraph  $\langle S_{l}\rangle$ is connected and, as a
consequence, the subgraph $\langle C_{S_l}\rangle$ also is
connected. Therefore, $|C_{S_l}|\leq |S_{l}|+1$.

Let $v\in C_{S_l}$. If  $S'=C_{S_l}\backslash \{v\}$ is a defensive
alliance in $\Gamma$, then $a(\Gamma)\le a({\cal L}(\Gamma))$.
Suppose  $S'=C_{S_l}\backslash \{v\}$ is not a defensive alliance in
$\Gamma$. In such case, there exists $u\in S'$ such that
$2|N_{S'}(u)|+2 \leq \delta (u).$ Since
$|N_{C_{S_l}}(u)|=N_{S'}(u)+1$, we have
\begin{equation}\label{Ss}
\delta (u)\ge 2|N_{C_{S_l}}(u)|.
\end{equation}
We shall use (\ref{Ss}) to show that  $S''=C_{S_l}\backslash \{u\}$
is a defensive alliance in $\Gamma$.

Suppose $w\in S''$ is a vertex adjacent to $u$ and let $e=\{u,w\}$.
Since $S_{l}$ is a defensive alliance in ${\cal L}(\Gamma)$,
$2|N_{S_l}(e)|+1\geq \delta(e)$. Therefore, by
$|N_{S_{l}}(e)|=|N_{C_{S_l}}(u)|+|N_{C_{S_l}}(w)|-2$ and $ \delta
(e)=\delta (u)+\delta(w)-2$, we obtain
\begin{equation}\label{Su}
2|N_{C_{S_l}}(u)|+2|N_{C_{S_l}}(w)|-1 \geq \delta (u)+\delta (w).
\end{equation}
By  (\ref{Ss}) and (\ref{Su}) we deduce $2|N_{C_{S_l}}(w)|- 1 \geq
\delta(w)$. Moreover, since   $|N_{C_{S_l}}(w)| =|N_{S''}(w)|+1$, we
have $2|N_{S''}(w)|+1\geq \delta(w).$ On the other hand, if  $w$ is
not adjacent to $u$, then $|N_{S''}(w)|=|N_{C_{S_l}}(w)|$. Hence,
$2|N_{S''}(w)|+ 1 \geq \delta(w)$. Thus, $S''$ is a defensive
alliance in $\Gamma$.
\end{proof}

It is easy to deduce sufficient conditions for $a(\Gamma)\leq
a({\cal L}(\Gamma))$ or $ a({\cal L}(\Gamma))\leq a(\Gamma)$ from
the above bounds and the bounds on $a(\Gamma)$ obtained in
\cite{alliancesOne,spectral}. For instance, it was shown in
\cite{alliancesOne} that $a(\Gamma)\le \left\lceil
\frac{n}{2}\right\rceil$. Hence, by Theorem \ref{thcota}, we have
$$\left\lceil\frac{n}{2}\right\rceil\le\left\lceil\frac{\delta_{n}+\delta_{n-1}-1}{2}\right\rceil
 \Rightarrow a(\Gamma)\le
{a}({\cal L}(\Gamma)).$$ In particular,
$$\frac{n}{2} < \delta_n \Rightarrow a(\Gamma)\le {a}({\cal L}(\Gamma)).$$

%The below result leads to a sufficient condition for $a(\Gamma)<
%a({\cal L}(\Gamma))$ and for  $ a({\cal L}(\Gamma))< a(\Gamma)$. The
%lower bound was obtained in \cite{spectral} and the upper bound was
%obtained in \cite{alliancesOne}.

%\begin{lemma} {\rm }\label{cotaconnectivity}
%Let $\Gamma$ be a simple graph of order $n$. Let $\mu$ be the
%algebraic connectivity of $\Gamma$. The defensive alliance number of
%$\Gamma$ is bounded by $$ \left\lceil\frac{n\mu}{n+\mu}\right\rceil
%\le a(\Gamma)\le \left\lceil \frac{n}{2}\right\rceil.$$
%\end{lemma}

%\begin{theorem} {\rm }
%Let $\Gamma$ be a simple graph of order $n$, size $m$ and maximum
%degree $\delta_1$. Let $\mu$ be the algebraic connectivity of
%$\Gamma$ and let $\mu_l$ be the algebraic connectivity of ${\cal
%L}(\Gamma$). Then,
%\begin{itemize}
%\item $\left\lceil
%\frac{n}{2}\right\rceil<\left\lceil\frac{m\mu_l}{m+\mu_l}\right\rceil
%\Rightarrow a(\Gamma)< a({\cal L}(\Gamma));$

%\item $\delta_1<\left\lceil\frac{n\mu}{n+\mu}\right\rceil \Rightarrow a({\cal L}(\Gamma))< a(\Gamma);$
%\end{itemize}
%\end{theorem}

%\begin{proof}
%The results are direct consequence of Lemma \ref{cotaconnectivity}
%and Theorem \ref{thcota}.
%\end{proof}
%Notice that analogous results can be derived in the case of strong
%alliances from $\left\lceil\frac{n(\mu+1)}{n+\mu}\right\rceil \le
%\hat{a}(\Gamma)\le \left\lceil \frac{n+1}{2}\right\rceil $ and
%Theorem \ref{thcota}.

\section{Global defensive alliances  and line graphs}

\begin{theorem}
Let $\Gamma$ be a simple graph of size $m>6$, then
$$\gamma_{a}({\cal L}(\Gamma))\ge \left\lceil
\sqrt{m+4}-1\right\rceil.$$
\end{theorem}

\begin{proof}
If  $S_{l}$ is a global defensive alliance in ${\cal L}(\Gamma)$,
then $$m-|S_{l}|\leq \sum_{v\in S_{l}}|N_{V_{l}\setminus
S_{l}}(v)|\leq \sum_{v\in
S_{l}}|N_{S_{l}}(v)|+|S_l|\leq|S_{l}|^{2}.$$
 On the other hand, if
 $|S_{l}|\le 2$, then
$|N_{V_l\setminus S_{l}}(v)|\le 2$, $\forall v\in S_{l}$. Thus,
$m\le 6$. Therefore, $ m>6\Rightarrow |S_{l}|>2$. By adding
$3\leq|S_{l}|$ and $m-|S_{l}|\leq |S_{l}|^{2}$, the result follows.
\end{proof}

The above bound is attained, for instance, in the case of the graph
of Figure \ref{figure}. In this case we can take the minimum global
defensive alliance as $S_l=\{a,b,g\}$.

 Several tight bounds on $\gamma_{a}({\cal
L}(\Gamma))$ and $\gamma_{\hat{a}}({\cal L}(\Gamma))$, in terms of
parameters of $\Gamma$, can be derived from the previous bounds on
$\gamma_{a}(\Gamma)$ and $\gamma_{\hat{a}}(\Gamma)$
\cite{planar,spectral,GlobalalliancesOne}. For instance, we consider
the following result.

\begin{theorem} {\rm \cite{spectral}}
Let $\Gamma$ be a simple graph of order $n$ and maximum degree
$\delta_1$. Then
$$\gamma_{{{a}}}(\Gamma)\ge \left\lceil \frac{2n}{\delta_1+3}\right\rceil \quad {\it and }\quad
\gamma_{{\hat{a}}}(\Gamma)\ge \left\lceil
\frac{n}{\left\lfloor\frac{\delta_1}{2}\right\rfloor+1}\right\rceil.$$
Both bounds are tight.
\end{theorem}

\begin{corollary}
Let $\Gamma$ be a simple graph of size $m$  whose maximum degrees
are  $\delta_{1}$ and $ \delta_{2}$. Then
$$\gamma_{a}({\cal L}(\Gamma))\geq \left\lceil\frac{2m}{\delta_{1}+\delta_{2}+1}\right\rceil
 \quad {\it and }\quad
 \gamma_{\hat{a}}({\cal L}(\Gamma))\geq \left\lceil\frac{2m}{\delta_{1}+\delta_{2}}\right\rceil .$$
\end{corollary}

In the case of connected alliances we obtain the following results.

\begin{theorem} \label{diametADG}
Let  $\Gamma=(V,E)$  be a connected graph of order $n$, size $m$ and
diameter $D(\Gamma)$. Then

 $ \gamma_{ca}(\Gamma)\ge
\left\lceil\sqrt{D(\Gamma)+n}-1 \right\rceil$ and $\gamma_{ca}
\left( {\cal L}(\Gamma)\right) \ge
\left\lceil\sqrt{D(\Gamma)+m-1}-1\right\rceil .$
\end{theorem}
\begin{proof}
If $S$ denotes a  global defensive alliance in $\Gamma$, then
$$
n-|S|\le \sum_{v\in S} | N_{V\setminus S}(v)| \le \sum_{v\in S} |
{N_S}(v) |+ |S|.
$$
On the other hand, if $S$ is a dominating set and  $\langle S
\rangle$ is connected, then $D(\Gamma)\leq D(\langle S \rangle)+2.$
Hence,
\begin{equation} \label{cotadiam}
D(\Gamma)\leq |S|+1.
\end{equation} By adding $n-|S|\leq |S|^{2}$ and
(\ref{cotadiam}) we obtain the bound on $\gamma_{ca}(\Gamma)$. The
bound on $\gamma_{ca} \left( {\cal L}(\Gamma)\right)$ follows from
the bound on  $\gamma_{ca}(\Gamma)$ and  $D(\Gamma)-1\le D({\cal
L}(\Gamma))$.
\end{proof}

Let $\Gamma$ be the left hand side graph of Figure \ref{figure}. The
set $S=\{1,2,3\}$ is a global defensive alliance in $\Gamma$ and
$\langle S\rangle$ is connected. On the other hand, $S_l=\{a,b,g\}$
is a global defensive alliance in ${\cal L}(\Gamma)$ and $\langle
S_l\rangle$ is connected. In this case, Theorem \ref{diametADG}
leads to $ \gamma_{ca}(\Gamma)\ge 3$ and $\gamma_{ca} \left( {\cal
L}(\Gamma)\right) \ge 3$. Thus, the bounds are tight.

\end{document}